\documentclass[a4paper,12pt]{article}
\usepackage{amssymb}
\usepackage{graphicx}
\def\beq{\begin{equation}}
\def\eeq{\end{equation}}
\def\bea{\begin{eqnarray}}
\def\eea{\end{eqnarray}}
\def\nn{\nonumber}

\def\braket#1#2{\VEV{#1 | #2}}
\def\VEV#1{\left\langle #1\right\rangle}
\def\tR{\tilde{R}}
\def\tmu{\tilde{\mu}}
\def\tS{\tilde{S}}
\newtheorem{theorem}{Theorem} 
\newtheorem{prop}{Proposition} 
\renewcommand{\labelenumi}{\rm (\roman{enumi})}

\begin{document}
\setlength{\baselineskip}{15pt}
\vspace*{2cm}
\begin{center}

\textsf{\Large All Link Invariants for Two Dimensional Solutions \\
of Yang-Baxter Equation and Dressings}

\bigskip
N. Aizawa\footnote{Present address : 
Department of Mathematics and Information Sciences, Osaka Prefecture University, 
Daisen Campus, Sakai, Osaka 590-0035, Japan.}, M. Harada, M. Kawaguchi and E. Otsuki \\

\bigskip
\textit{
Department of Applied Mathematics,\\
Osaka Women's University, \\
Sakai, Osaka 590-0035, Japan
}
\end{center}

\vfill
\begin{abstract}
All polynomial invariants of links for two dimensional solutions of 
Yang-Baxter equation is constructed by employing Turaev's method. 
As a consequence, it is proved that the best invariant so constructed is 
the Jones polynomial and there exist three solutions connecting to 
the Alexander polynomial. Invariants for higher dimensional solutions, 
obtained by the so-called dressings, are also investigated. 
It is observed that the dressings 
do not improve link invariant unless some restrictions are put on dressed solutions.
\end{abstract}
\clearpage
%
%
%
%
\section{Introduction}	

  Importance of Yang-Baxter equation (YBE) in knot theory has been widely recognized.  
In the context of knot theory, the YBE appears as a matrix expression of Reidemeister move 
of type III when we assign a matrix to each crossing of a link diagram. 
This means that the YBE is an essential element for constructing topological invariants of 
links. Turaev gave a general method to construct link invariants starting from 
solutions of YBE \cite{Tu}. Kauffman gave a more diagrammatic method to construct 
invariants by interpreting a solution of YBE as vertex weights of state expansion of 
a given link diagram \cite{LK}. Solutions of YBE also appear in statistical mechanics 
as a Boltzmann weight of exactly solvable models. 
Motivated by this fact, Akutsu, Wadati and Deguchi developed a general method of 
constructing link invariants starting from exactly solvable models in 
statistical mechanics [3-9]. 

  In any aforementioned construction of invariants, we need to know solutions 
of YBE explicitly. One way of finding solutions of YBE is to use representations 
of quantum groups (see for example \cite{CP}). In general, for each representation 
of a quantum group, there exist an associated solution of YBE. 
This method, therefore, can give arbitrary dimensional solutions. 
However this construction does not exhaust all possible solutions for fixed dimension. 
Classification of solutions of YBE has not been completed yet except constant case 
in two dimension. We mean by \textit{constant} a YBE and its solutions which do not 
contain spectral parameters. 
Two dimensional constant solutions of YBE have been classified by Hietarinta \cite{Hi1,Hi2}. 
While classification of solutions having spectral parameters was made in \cite{SUAW} where 
eight-vertex type solutions of YBE in two dimension were studied in connection with 
solvable two-component models. 

In the present work, we investigate all possible link invariants obtained 
from the two dimensional constant solutions of YBE. We apply Turaev's construction to 
the two dimensional solutions obtained by Hietarinta. We anticipated to find new 
polynomial invariants of links, since two dimensional solutions contain at most three 
independent parameters. However the present investigation reveals that 
two dimensional solutions do not produce better invariants than already known ones. 
It is also shown that three two dimensional $R$-matrices are related to the Alexander polynomial. 

  The above mentioned results motivate us to study higher dimensional constant solutions. 
We study higher dimensional solutions obtained not via quantum groups 
but via the method called $dressing$ \cite{Hi3}. 
Dressing is a method to construct higher dimensional solutions 
from known lower dimensional ones.  We shall show that dressing do not improve 
link invariants unless we put some additional conditions on dressing. 

  This paper is organised as follows. We give a brief review of Turaev's construction 
and two dimensional solutions of YBE in the next two sections in order to fix our notations 
and conventions. Section 4 is devoted to construction and classification of link invariants 
for two dimensional solutions. We will see that there exist three $R$-matrices 
from which the Alexander polynomial is obtained. A state model for one of them is 
discussed also in section 4. Link invariants for higher dimensional solutions 
obtained by dressings are investigated in section 5. Some concluding remarks 
are given in section 6.

%
%
%
\section{Turaev's Construction of Link Invariants}	

  Let $V$ be a $N$ dimensional vector space over a filed $\mathbb K,$ 
and $ R :  V \otimes V  \ \rightarrow \ V \otimes V $ 
satisfy the YBE of the braid group form
\beq
  R_{12} R_{23} R_{12} = R_{23} R_{12} R_{23}. \label{YBE}
\eeq
An enhanced Yang-Baxter operator (EYB) is a collection 
$( R, \mu, \alpha, \beta),$ where $ \mu : V \ \rightarrow \ V $ 
is a homomorphism and $ \alpha, \beta $ are invertible elements 
of $\mathbb K,$ subject to the conditions 
\bea
 & & R \circ (\mu \otimes \mu) = (\mu \otimes \mu) \circ R, \label{EYB1} \\
 & & {\rm Tr}_2(R \circ (\mu \otimes \mu)) = \alpha \beta \mu, \label{EYB2} \\
 & & {\rm Tr}_2(R^{-1} \circ (\mu \otimes \mu)) = \alpha^{-1} \beta \mu, \label{EYB3}
\eea 
where ${\rm Tr}_2$ is a trace for the second vector space in the tensor product. 
It is widely known that the $R$-matrix gives rise to a representation of 
the $n$-strand braid group $B_n$ on $ V^{\otimes n}.$ 
Let $ \{ \ \sigma_i,\ i=1,\cdots,n-1 \ \} $ be generators of $ B_n. $ 
Define an isomorphism $ R_i : V^{\otimes n} \ \rightarrow \ V^{\otimes n} $ by 
\beq
  R_i = {\rm id}^{\otimes(i-1)} \otimes R \otimes {\rm id}^{\otimes(n-i-1)}. \label{Ri}
\eeq
Then the representation $(\pi,V^{\otimes n})$ of $ B_n$ is defined by 
$ \pi(\sigma_i^{\pm 1}) = R_i^{\pm 1}. $ 
Using the representations, for each EYB $ S = (R, \mu, \alpha, \beta) $ of  
a given $R$-matrix, an invariant of links  immediately follows. 
Let $L$ be an oriented link, $ \xi $ be a $n$-strand braid whose closure $\bar{\xi}$ is 
isotopic to $L.$  One can show that the mapping 
$ {\displaystyle T_S : \bigcup_{n \geq 1} B_n \ \rightarrow \ {\mathbb K} }$ defined 
for the braids gives rise to an invariants of links
\beq
   T_S(L) = T_S(\xi) = \alpha^{-w(\xi)} \beta^{-n} {\rm Tr}(\pi(\xi) \circ \mu^{\otimes n}),
   \label{Ts}
\eeq 
where $w(\xi)$ is the writhe of $ \xi $ defined by $ w(\xi) =  w(\bar{\xi}). $  

\begin{figure}[htbp] 
  \centerline{
  \includegraphics[scale=0.6, trim=0 0 0 0]{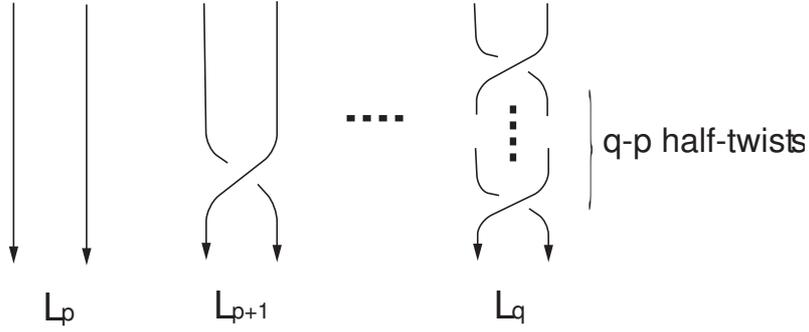}
  } 
  \caption{Links $ L_i.$ The diagrams are assumed to be 
  identical excepting the depicted parts.}
  \label{SkeinD}
\end{figure}

  The invariant $T_S(L)$ has the following properties
\begin{enumerate}
  \item If a link $L$ is the disjoint union of two links $L_1 $ and $ L_2$, 
    then
    \beq
       T_S(L) = T_S(L_1) T_S(L_2). \label{TSL1L2} 
    \eeq
  \item If the $R$-matrix satisfies the relation
    \beq
       \sum_{i=p}^q k_i R^i = 0, \qquad k_p, \cdots, k_q \in {\mathbb K}
       \label{annihC}
    \eeq
    then the invariants $T_S$ satisfies the skein-type relation
    \beq
       \sum_{i=p}^q k_i \alpha^i T_S(L_i) = 0,   \label{skein}
    \eeq
    where $L_i$'s are links depicted in Figure~\ref{SkeinD}.
  \item If $ S = (R, \mu, \alpha, \beta) $ is an EYB, then  
    $ S_1 = (-R,-\mu,\alpha,\beta),\ S_2 = (R,\mu,-\alpha,-\beta) $ 
    and $ S_3 = (R,-\mu,-\alpha,\beta) $ are also EYB.   
    For any $\ell$-component link $L$, the corresponding invariants 
    are related
    \beq
      T_{S_1}(L) = T_{S_2}(L) = T_{S_3}(L) = (-1)^{\ell} T_S(L).  \label{TS123}
    \eeq
\end{enumerate}

%
%
%
\section{Two Dimensional Solutions of YBE}	
\setcounter{equation}{0}

Two dimensional ($ N=2$) constant solutions to YBE are extensively studied by Hietarinta 
by using computer \cite{Hi1,Hi2}. 
The classification of the solutions was made up to the invariance 
of YBE by the following transformations
\bea
  & & R \quad\ \ \, \rightarrow \quad R' = \kappa (Q \otimes Q) R ( Q \otimes Q )^{-1},
  \label{similarity} \\
  & & R^{k\ell}_{ij} \quad \rightarrow \quad R'^{k\ell}_{ij} = R^{ij}_{k\ell}, 
  \label{transpose} \\
  & & R^{k\ell}_{ij} \quad \rightarrow \quad R'^{k\ell}_{ij} = R^{k+n, \ell+n}_{i+n, j+n} \ \ 
      {\rm (indices\ mod }\ N) ,
  \label{shift} \\
  & & R^{k\ell}_{ij} \quad \rightarrow \quad R'^{k\ell}_{ij} = R^{\ell k}_{ji},
  \label{LRflip}
\eea
where $ \kappa \neq 0 \in {\mathbb K} $ and $ Q $ is a nonsingular $ N \times N $ matrix. 
These transformations do not affect the invariant $T_S(L)$ defined in the last section. The proof of 
this fact is given in Appendix~A. It is, therefore, legitimate to investigate link invariants 
according to the classification of \cite{Hi1,Hi2}. 

  Hietarinta obtained 35 solutions (including zero matrix); one three-parameter, 
four two-parameter, 15 one-parameter and 14 no-parameter solutions.
No-parameter solutions will give numerical invariants which are less powerful 
than polynomial invariants. Thus we exclude such solutions. The solutions given 
by singular matrices are also excluded, since inverse matrices are assigned to 
negative crossings. We give a list of nonsingular solutions excepting no-parameter 
ones. Entries of the $R$-matrices are labelled as follows
\beq
  R = \left(
  \begin{array}{cccc}
    R_{++}^{++} & R_{++}^{-+} & R_{++}^{+-} & R_{++}^{--} \\
    R_{-+}^{++} & R_{-+}^{-+} & R_{-+}^{+-} & R_{-+}^{--} \\
    R_{+-}^{++} & R_{+-}^{-+} & R_{+-}^{+-} & R_{-+}^{--} \\
    R_{--}^{++} & R_{--}^{-+} & R_{--}^{+-} & R_{--}^{--}
  \end{array}
  \right).
  \label{Rsuff}
\eeq

\medskip \noindent
3 parameters
$$
  R_{3.1} = \left(
  \begin{array}{cccc}
     1 & 0 & 0 & 0 \\
     0 & 0 & q & 0 \\
     0 & p & 0 & 0 \\
     0 & 0 & 0 & s
  \end{array}
  \right).
$$

\medskip\noindent
2 parameters
$$
     R_{2.1} = \left(
     \begin{array}{cccc}
        1 & 0 & 0 & 0 \\
        0 & 0 & q & 0 \\
        0 & p & 1-pq & 0 \\
        0 & 0 & 0 & 1
     \end{array}
     \right),
     \qquad
          R_{2.2} = \left(
     \begin{array}{cccc}
        1 & 0 & 0 & 0 \\
        0 & 0 & q & 0 \\
        0 & p & 1-pq & 0 \\
        0 & 0 & 0 & -pq
     \end{array}
     \right),
$$
$$
     R_{2.3} = \left(
     \begin{array}{cccc}
        1 & 1 & p & q \\
        0 & 0 & 1 & 1 \\
        0 & 1 & 0 & p \\
        0 & 0 & 0 & 1
     \end{array}
     \right).
$$

\medskip\noindent
1 parameter
$$
     R_{1.1} = \left(
     \begin{array}{cccc}
        1+2q-q^2 & 0 & 0 & 1-q^2 \\
        0 & 1-q^2 & 1+q^2 & 0 \\
        0 & 1+q^2 & 1-q^2 & 0 \\
        1-q^2 & 0 & 0 & 1-2q-q^2
     \end{array}
     \right),
$$
$$
     R_{1.2} = \left(
     \begin{array}{cccc}
        1 & 0 & 0 & 1 \\
        0 & 0 & q & 0 \\
        0 & 1 & 1-q & 0 \\
        0 & 0 & 0 & -q
     \end{array}
     \right),
     \qquad
          R_{1.3} = \left(
     \begin{array}{cccc}
        1 & 1 & -1 & q \\
        0 & 0 & 1 & -q \\
        0 & 1 & 0 & q \\
        0 & 0 & 0 & 1
     \end{array}
     \right),
$$
$$
     R_{1.4} = \left(
     \begin{array}{cccc}
        0 & 0 & 0 & q \\
        0 & 1 & 0 & 0 \\
        0 & 0 & 1 & 0 \\
        q & 0 & 0 & 0
     \end{array}
     \right).
$$

%
%
%
\section{Link Invariants}	
\setcounter{equation}{0}
\subsection{EYB and List of Invariants}

  In this section, we derive all invariants of links obtained from the $R$-matrices 
listed in the previous section. Key of the derivation is to find an EYB for a given 
$R$-matrix. Our procedure to obtain all possible EYB for a given $R$-matrix 
is summarised as follows.
\begin{enumerate} 
  \item Regarding the condition (\ref{EYB2}) as a set of equations with  
  $ \mu, \alpha$ and $ \beta $  as unknown variables, 
  solve it with MAPLE. Usually, this gives several solutions containing 
  some variables that are not determined yet.
  \item For each solution, the undetermined variables are fixed so as to satisfy the 
  conditions (\ref{EYB1}) and (\ref{EYB3}). There are some cases where the conditions  
  require unacceptable results; $ \mu = 0, $ or $ \alpha = 0, $ or $ \beta = 0, $ or 
  all the parameters in 
  the $R$-matrix being fixed. These cases are excluded. 
\end{enumerate}
The EYB's so obtained are summarised in Table~\ref{EYBtable}. 
Note that one $R$-matrix gives rise to some inequivalent EYB's. 

\begin{table}[htbp]
\caption{EYB and link invariants. $\beta$ takes arbitrary value for all cases.}
\medskip
\centerline{
\begin{tabular}{|c | c | c | c | c|}
\hline
 $R$ & $\mu$ & $\alpha$ & invariants & comments \\
\hline
 {}  & diag$(\pm \beta, \pm \beta)$ & $ \pm 1$ & $1$ (for knots) & $s=1$ \\
\cline{2-5}
 $R_{3.1}$ & diag$(\pm \beta,\mp \beta)$ & $ \pm 1$ & $0$ (for knots) & $s=-1$ \\
\cline{2-5}
 {} & ${\displaystyle 
              \begin{array}{c}
                {\rm diag}(\pm \beta, 0)
               \\
                {\rm diag}(0, \pm \beta)
              \end{array}} $  
            & $ {\displaystyle 
                \begin{array}{c}
                 \pm 1 \\
                 \pm s 
                \end{array}
                }$ 
             & \begin{tabular}{cc}
                 $1$ \\ $1$
               \end{tabular}
             & \\
\hline
 {}          & diag$(\pm\sqrt{pq}\beta, \pm\beta/\sqrt{pq}) $ & $\pm\frac{1}{\sqrt{pq}}$  & Jones & \\
 {}          & diag$(\pm \beta, 0)$ & $\pm 1$ & $1$ & \\
 {}          & diag$(0, \pm \beta)$ & $\pm 1$ & $1$ & \\
 \cline{2-5}
 $ 
 \begin{array}{c}
 R_{2.1} \\ \\
 \end{array}
 $           
              & ${\displaystyle 
                \left(\begin{array}{cc}
                      \pm \beta & 0 \\ \lambda & 0
                \end{array}\right)
               }$
              & $\pm 1$ 
              & $1$
              & $\begin{array}{c} \\ q=1 \end{array}$ \\
  {}          & $ {\displaystyle
                   \left(\begin{array}{cc}
                      0 & \lambda \\ 0 & \pm \beta
                   \end{array}\right)
                }$
              & $\pm 1$ 
              & $1$
              & \begin{tabular}{c}$\lambda$ is arbitrary \\ \\ \end{tabular} \\
\hline
 $R_{2.2}$   &$ {\displaystyle
                \begin{array}{c}
                    {\rm diag}(\pm\beta/\sqrt{pq}, \mp\beta/\sqrt{pq})
                  \\
                    {\rm diag}(\pm \beta, 0)
                  \\
                    {\rm diag}(0, \mp \beta)
                \end{array}
             }$
             & $ {\displaystyle
                \begin{array}{c}
                  \pm\sqrt{pq}
                  \\
                  \pm 1
                  \\
                  \pm pq
                \end{array}
             }$
             & \begin{tabular}{c}
                 $0$ \\ $1$ \\ $1$
               \end{tabular}
             & \begin{tabular}{c}
                 Alexander \\ \\ \\
               \end{tabular}
             \\
\hline
 $R_{2.3}$   & diag$(\pm \beta,\pm \beta)$ & $\pm 1$ & $ 2^{\ell}$ & $p=-1$ \\
\hline
 {}          & diag$(\pm \beta,\mp\beta)$ & $ \pm 2q $ & $0$ & Alexander \\
 {}   & $ {\displaystyle
                  \pm \frac{\beta}{2}
                  \left(\begin{array}{cc}
                    1+q & \sqrt{1-q^2} \\ \sqrt{1-q^2} & 1-q
                  \end{array}\right)
               }$
             & $ \pm 2 $ & $1$ & \\
  $R_{1.1}$  & $ {\displaystyle  
                  \pm \frac{\beta}{2}
                  \left(\begin{array}{cc}
                    -1-q & \sqrt{1-q^2} \\ \sqrt{1-q^2} & -1+q
                  \end{array}\right)
              }$
             & $ \mp 2$ & $1$  & \\
  {}         & $ {\displaystyle 
                  \pm \frac{\beta}{2q}
                  \left(\begin{array}{cc}
                    1+q & \sqrt{1-q^2} \\ -\sqrt{1-q^2} & -1+q
                  \end{array}\right)
               } $
             & $ \pm 2 $ & $1$  & \\
  {}         & $ {\displaystyle 
                  \pm \frac{\beta}{2q}
                  \left(\begin{array}{cc}
                    -1-q & \sqrt{1-q^2} \\ -\sqrt{1-q^2} & 1-q
                  \end{array}\right)
               }$
             & $ \mp 2$ & $1$ & \\
\hline
  {}         & diag$(\pm\beta/\sqrt{q}, \mp\beta/\sqrt{q}) $ 
             & $ \pm \sqrt{q} $ & $0$ & Alexander \\
  $R_{1.2}$  & $ {\displaystyle
               \left(\begin{array}{cc}
                  \pm \beta & \pm \beta/\sqrt{1+q} \\ 0 & 0
               \end{array}\right)
               }$
             & $ \pm 1 $ & $1$  & \\
  {}         & $ {\displaystyle
               \left(\begin{array}{cc}
                  \pm \beta & \mp \beta/\sqrt{1+q} \\ 0 & 0 
               \end{array}\right)
               } $
             & $ \pm 1 $ & $1$  & \\
\hline
  $ R_{1.3}$ & $ {\displaystyle
               \left(\begin{array}{cc}
                  \pm \beta & \mp (1+q) \beta \\ 0 & \pm \beta
               \end{array}\right)
               } $
             & $ \pm 1 $ & $2^{\ell}$ & \\
\hline
 $ R_{1.4}$  & diag$(\pm\beta,\pm\beta)$ & $ \pm 1$ & $1$ (for knots) & \\
\hline
\end{tabular}}
\label{EYBtable}
\end{table}

  In Table~\ref{EYBtable}, the values of link invariants corresponding to 
each EYB are also shown. Each EYB has two possibilities of signs. However, 
because of (\ref{TS123}), the sign difference of EYB causes a difference 
of overall factor of invariants. The indicated invariants correspond to 
one of the sign choices. Some details on the invariants are studied in the 
following.  
\begin{enumerate}
\item $ R_{3.1} $ with $ s=1,$ $ \mu = {\rm diag}(\beta,\beta),\ \alpha = 1. $
 
   For this choice of $\mu $ and $ \alpha,$ $ s=1$ is required to 
  satisfy the conditions of EYB. Thus the $R$-matrix $ R = R_{3.1} $ 
  is reduced to two-parameter and it satisfies
  \beq
     R^2 - R -pq  + pq R^{-1} = 0. \label{R31rel}
  \eeq
  By (\ref{skein}), the invariant satisfies the skein-type relation
  \beq
     T_S(L_{++}) - T_S(L_+) - pq T_S(L_0) + pqT_S(L_-) = 0,
     \label{R31skein}
  \eeq
  where $ L_{++}, L_+, L_0 $ and $ L_- $ are links depicted in 
  Figure~\ref{Lpp}.
  \begin{figure}[htbp]
    \begin{center}
      \includegraphics[scale=0.6]{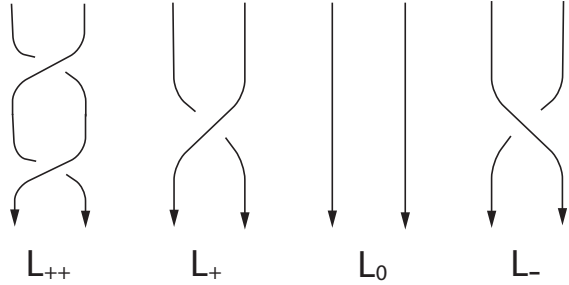}
    \end{center}
    \caption{Links $ L_{++}, L_+, L_0$  and $ L_-.$}
    \label{Lpp}
  \end{figure}
  The value for the unknot $U$ is $ T_S(U) = 2.$ Furthermore, it turns 
  out that $ T_S(K) = 2 $ for any knot $K$.  To see this, we derive 
  another version of skein-type relation.
  Multiplying $R$ to (\ref{R31rel}) and doing a sum of it and (\ref{R31rel}), 
  we obtain
  \[
    R^3 - (1+pq)R + pqR^{-1} = 0.
  \]
  Thus the invariant satisfies
  \beq
     T_S(L_{+++}) = (1+pq) T_S{L_+} - pq T_S(L_-). \label{R31skeinal} 
  \eeq
  The induction on the number of crossings proves that $ T_S(K) = 2 $ for 
  any knot $K.$ By normalizing $T_S(K) $ by the value for unknot, 
  we obtain $ T_S(K) = 1.$ 
  \item $ R_{3.1} $ with $ s = -1,\ \mu = {\rm diag}(\beta,-\beta),\ \alpha = 1. $ 
  
     For this choice of $\mu $ and $ \alpha,$ $ s=-1$ is required to 
  satisfy the conditions of EYB. Thus the $R$-matrix $ R = R_{3.1} $ 
  is reduced to two-parameter and it satisfies
  \beq
    R^2 - (1+pq) + pq R^{-2} = 0. \label{R31rel2}
  \eeq
  Thus we obtain
  \beq
     T_S(L_{++}) - (1+pq)T_S(L_0) + pq T_S(L_{--}) = 0. \label{R31skein2}
  \eeq
  It turns out that $ T_S(K) = 0 $ for any knot $K.$ One can prove it in a 
  way similar to case (i). The reason of the vanishing invariants is that  
  the invariant for unknot is equal to 0. 
  Even if we regularize the value for unknot being nonvanishing, one can show, by the 
  skein relation, that $T_S$ for trefoil is the same as the one for unknot. Thus the invariant  
  is less powerful than the Jones or the Alexander polynomial.
  \item $ R_{3.1},\ \mu = {\rm diag}(\beta,0),\ \alpha = 1. $ 
  
    There is no restriction on $R = R_{3.1} $ in this case. It is easy to verify 
  \beq
     R^2 - (1+s)R + (s-pq) + pq(1+s)R^{-1}-spqR^{-2} = 0. \label{R31rel3}
  \eeq
  It follows that
  \beq
     T_S(L_{++}) - (1+s) T_S(L_+) + (s-pq) T_S(L_0) 
     + pq(1+s) T_S(L_-) -spq T_S(L_{--}) = 0. \label{R31skein3}
  \eeq
  Although the skein-type relation has a complex form, it turns out that 
  $ T_S(L) = 1 $ for any link $L.$ This is a consequence of 
  $ (\mu \otimes \mu) R^{\pm 1} = \mu \otimes \mu $ and the definition of 
  invariant (\ref{Ts}). 
  \item $ R_{3.1},\ \mu = {\rm diag}(0,\beta),\ \alpha = s. $
  
    This case is exactly the same as case (iii). 
    The skein-type relation is  given by (\ref{R31skein3}) and 
    $ T_S(L) = 1 $ for any link $L,$ 
    since $ (\mu \otimes \mu) R^{\pm 1} = s^{\pm 1} (\mu \otimes \mu). $ 
  \item $ R_{2.1},\ \mu = {\rm diag}(\sqrt{pq}\beta,\beta/\sqrt{pq}), \ \alpha = 1/\sqrt{pq}. $
  
     The matrix $ R = R_{2.1}$ satisfies the relation
    \beq
        R + (pq-1) - pq R^{-1} = 0. \label{R21rel}
    \eeq 
    It follows that
    \beq
       (pq)^{-1} T_S(L_+) - pq\, T_S(L_-) = (1/\sqrt{pq} - \sqrt{pq})\, T_S(L_0). \label{R21skein}
    \eeq 
    The value for the unknot is $ T_S(U) = \sqrt{pq}+1/\sqrt{pq}. $ 
    The two parameters in $R$ are combined together to give a polynomial in one variable $t=pq.$ 
    This invariant is the Jones polynomial. 
   \item $ R_{2.1}, \ \mu = {\rm diag}(\beta,0),\ \alpha = 1. $
   
     The skein relation is found from (\ref{R21rel}) and it is given by
     \beq
        \frac{1}{\sqrt{pq}} T_S(L_+) - \sqrt{pq}\, T_S(L_-) = (\frac{1}{\sqrt{pq}} - \sqrt{pq})\, 
        T_S(L_0).  \label{R21skein2}
     \eeq
     Although the skein relation has the form of specialisation of HOMFLY polynomial, 
     the invariant takes the trivial value; $ T_S(L) = 1 $ for any link $L$. This is due to 
     the same reason as case (iii). 
    \item $ R_{2.1}, \ \mu = {\rm diag}(0,\beta),\ \alpha = 1. $ 
    
      This case is exactly the same as case (vi).
    \item $ R_{2.1} $ with $ q = 1,\ {\displaystyle \mu = \left(
           \begin{array}{cc}
             \beta & 0 \\ \lambda & 0
           \end{array}
           \right),\ \alpha = 1.} $
        
       For this choice of $\mu $ and $ \alpha,$ $ q=1$ is required to 
     satisfy the conditions of EYB. Thus the $R$-matrix $ R = R_{2.1} $ 
     is reduced to one-parameter. The matrix $\mu$ contains an additional 
     free parameter $ \lambda. $ The skein relation is obtained by setting 
     $ q = 1 $ in (\ref{R21skein2}). It turns out again that $ T_S(L) = 1 $ 
     for any link $L.$ This is due to the same reason as case (iii).
    \item $ R_{2.1} $ with $ q = 1,\ {\displaystyle \mu = \left(
           \begin{array}{cc}
             0 & \lambda \\ 0 & \beta
           \end{array}
           \right),\ \alpha = 1.} $
       
       This case is exactly the same as case (viii).
    \item $ R_{2.2}, \ \mu = {\rm diag}(\beta/\sqrt{pq},-\beta/\sqrt{pq}),\ \alpha = \sqrt{pq}.$ 
    
       The $R=R_{2.2}$ satisfies the same relation as $R_{2.1}$
       \beq
          R + (pq-1) - pq R^{-1} = 0. \label{R22rel}
       \eeq 
       It follows that the skein relation is 
       \beq
          T_S(L_+) - T_S(L_-) = (\frac{1}{\sqrt{pq}}-\sqrt{pq}) T_S(L_0). \label{R22skein}
       \eeq
       This is the skein relation of Alexander polynomial. However, it turns out that 
       $ T_S(L) = 0 $ for any link $L.$ 
       To see this, note that if a link $L$ is the disjoint union of two links $L_1$ and $L_2$, 
       then, because of the skein relation (\ref{R22skein}), $ T_S(L) = 0. $ While by (\ref{TSL1L2}) 
       $ T_S(L) = T_S(L_1) T_S(L2) = 0. $ This means $ T_S(L_1) = 0 $ or $ T_S(L_2) = 0. $ 
       Thus $ T_S(L) = 0 $ for any $L.$ 
    \item $ R_{2.2}, \ \mu = {\rm diag}(\beta,0),\ \alpha = 1.$
    
       Because of the same reason as case (iii), $ T_S(L) = 1 $ for any link $L.$
       
    \item $ R_{2.2}, \ \mu = {\rm diag}(0,\beta),\ \alpha = -pq.$

       Note that $ (\mu \otimes \mu) R^{\pm 1} = -(pq)^{\pm 1} (\mu \otimes \mu) $ 
       holds in this case. 
       The same argument as case (iv) shows that $ T_S(L) = 1 $ for any link $L.$ 
    \item $ R_{2.3}$ with $ p=-1$, $\mu = {\rm diag}(\beta,\beta),\ \alpha = 1.$  
    
       For this choice of $\mu $ and $ \alpha,$ $ p=-1$ is required to 
      satisfy the conditions of EYB. Thus the $R$-matrix $ R = R_{2.3} $ 
      is reduced to one-parameter and it satisfies
      \beq
        R^2 - R - 1  + R^{-1} = 0. \label{R23rel}
      \eeq
      It follows that the invariant satisfies the skein-type relation
      \beq
         T_S(L_{++}) - T_S(L_+) - T_S(L_0) + T_S(L_-) = 0. \label{R23skein}
      \eeq
      One can see,  by the induction on the number of crossings, 
      that $ T_S(L) = 2^{\ell} $ for a $\ell$-component link $L.$ 
    \item $ R_{1.1}, \ \mu = {\rm diag}(\beta,-\beta), \ \alpha = 2q. $
      The matrix $R = R_{1.1}$ satisfies the relation
      \beq
         R + 2(q^2-1) - 4q^2 R^{-1} = 0. \label{R11rel}
      \eeq
      The skein relation is given by
      \beq
         T_S(L_+) - T_S(L_-) = (q^{-1}-q) T_S(L_0).  \label{R11skein}
      \eeq
      This is the skein relation of the Alexander polynomial. 
      However, by the same reason as case (x), it turns out that 
      $ T_S(L) = 0 $ for any link $L.$ 
    \item $ R_{1.1}, \ {\displaystyle \mu = \frac{\beta}{2}
           \left(
             \begin{array}{cc}
                1+q & \sqrt{1-q^2} \\ \sqrt{1-q^2} & 1-q
             \end{array}
           \right), \alpha = 2.} $
           
       The skein relation follows from (\ref{R11rel}) and it reads
       \beq
           q^{-1} T_S(L_+) - qT_S(L_-) = (q^{-1}-q) T_S(L_0). \label{R11skein2}
       \eeq
       It is easy to see that 
       $ (\mu \otimes \mu) R^{\pm 1} = 2^{\pm 1} (\mu \otimes \mu). $ 
       Therefore, by the same reason as case (iv), the invariant takes 
       a single value $T_S(L) = 1 $ for any link $L.$ 
    \item $ R_{1.1}, \ {\displaystyle \mu = \frac{\beta}{2}
           \left(
             \begin{array}{cc}
                1+q & -\sqrt{1-q^2} \\ -\sqrt{1-q^2} & 1-q
             \end{array}
           \right), \alpha = 2.} $

        This case is exactly the same as case (xv).
    \item $ R_{1.1}, \ {\displaystyle \mu = \frac{\beta}{2q}
           \left(
             \begin{array}{cc}
                1+q & \sqrt{1-q^2} \\ -\sqrt{1-q^2} & -1+q
             \end{array}
           \right), \alpha = 2.} $

       This case is exactly the same as case (xv).
    \item $ R_{1.1}, \ {\displaystyle \mu = \frac{\beta}{2q}
           \left(
             \begin{array}{cc}
                1+q & -\sqrt{1-q^2} \\ \sqrt{1-q^2} & -1+q
             \end{array}
           \right), \alpha = 2.} $

       This case is exactly the same as case (xv).
    \item $ R_{1.2},\ \mu = {\rm diag}(\beta/\sqrt{q},-\beta/\sqrt{q}),\ \alpha = \sqrt{q}. $
    
      The matrix $ R = R_{1.2}$ satisfies the relation
      \beq
        R + (q-1) - qR^{-1} = 0. \label{R12rel}
      \eeq
      It follows that the invariant satisfies the skein relation
      \beq
        T_S(L_+) - T_S(L_-) = (\frac{1}{\sqrt{q}} - \sqrt{q}) T_S(L_0). \label{R12skein}
      \eeq
      This is the skein relation of the Alexander polynomial. 
      However, by the same reason as case (x), it turns out that 
      $ T_S(L) = 0 $ for any link $L.$ 
    \item $ R_{1.2}, \ {\displaystyle \mu =            
           \left(
             \begin{array}{cc}
                \beta & \beta/\sqrt{1+q} \\ 0 & 0
             \end{array}
           \right), \alpha = 1.} $
           
       The skein relation follows from (\ref{R12rel}) and it reads
       \beq
          \frac{1}{\sqrt{q}} T_S(L_+) - \sqrt{q} T_S(L_-) = (\frac{1}{\sqrt{q}}-\sqrt{q}) T_S(L_0).
          \label{R12skein2}
       \eeq
       It turns out that $ T_S(L) = 1 $ for any link $L,$ due to the same reason as case (iii).
    \item $ R_{1.2}, \ {\displaystyle \mu =            
           \left(
             \begin{array}{cc}
                \beta & -\beta/\sqrt{1+q} \\ 0 & 0
             \end{array}
           \right), \alpha = 1.} $

        This case is exactly the same as case (xx).
    \item $ R_{1.3}, \ {\displaystyle \mu =            
           \left(
             \begin{array}{cc}
                \beta & -(1+q) \beta \\ 0 & \beta
             \end{array}
           \right), \alpha = 1.} $

    The matrix $R = R_{1.3} $ is peculiar, since it satisfies $ R^2 = 1. $ 
    It follows that the invariant has the property $ T_S(L_+) = T_S(L_-). $ 
    The value for the unknot is $ T_S(U) = 2. $ 
    Thus  $ T_S(L) = 2^{\ell} $ for a $\ell$-component link $L$.
    \item $ R_{1.4}, \ \mu = {\rm diag}(\beta,\beta),\ \alpha = 1. $

      The matrix $ R = R_{1.4} $ satisfies the relation
     \beq
        R^2 - R - q^2 + q^2 R^{-1} = 0. \label{R14rel}
     \eeq
     It follows that the skein-type relation is given by
     \beq
        T_S(L_{++}) - T_S(L_+) - q^2 T_S(L_0) + q^2 T_S(L_-) = 0. 
        \label{R14skein}
     \eeq
     One can show, in a  way similar to case (i), that $ T_S(K) = 1 $ 
     for any knot $K.$ 
\end{enumerate}
Besides the above listed $R$-matrices, there exist four non-singular 
no-parameter $R$-matrices. It is seen that one can construct EYB's from 
those four $R$-matrices and they give numerical invariants. 
We do not show them here, because numerical invariants are less 
powerful than polynomial ones. In summary, we obtained the 
following theorem.

\medskip
\begin{theorem}
 \label{EYBthm}\textit{
   All the non-singular two dimensional solutions of YBE gives rise to 
  EYB's. The best invariant of links obtained from the EYB's is the 
  Jones polynomial. 
  }
\end{theorem}

\subsection{A Yang-Baxter State Model for Alexander Polynomial}

  From Table~1, it is seen that three $R$-matrices, 
$ R_{2.2},\ R_{1.1} $ and $ R_{1.2}, $ are related to the Alexander polynomial. 
However, the Turaev's construction gives vanishing value for the corresponding 
invariants. We need to use another method to show that those $R$-matrices 
can be a source of the Alexander polynomial. The case of one-parameter 
restriction of $ R_{2.2} $ has been studied in some works. 
With a $R$-matrix obtained from $R_{2.2}$ by a scaling, 
Jaeger constructed a state model for the Alexander polynomial 
which can be interpreted as an ice-type model \cite{Ja}. 
This model was reinterpreted as a Yang-Baxter state model by 
Kauffman \cite{LK}. The supersymmetric aspects of the Alexander 
polynomial were studied in connection with the supersymmetric 
counterpart of a one-parameter restriction of $R_{2.2}$ \cite{KS,RS}. 
In \cite{KS}, the Alexander polynomial is described as a fermionic 
integral, that is, a state model involving bosonic and fermionic loops. 
The quantum field theory descriptions of the multi-variable 
Alexander polynomial based on the WZW model and the Chern-Simons model 
are given in \cite{RS}. 

  These works encourage us to study $ R_{1.1} $ and $ R_{1.2} $ in the 
same footing. We note that there is a fundamental difference between 
$ R_{1.1}, R_{1.2} $ and $ R_{2.2} $ besides the number of parameters. 
The matrix $R_{2.2}$ is \textit{spin preserving}, however $R_{1.1}$ and 
$R_{1.2}$ are not. The $R$-matrices having the following property 
are said to be spin preserving 
\beq
  R^{ab}_{cd} \neq 0, \quad {\rm if \ and \ only \ if} \quad a+b=c+d.
  \label{spinpre}
\eeq
Since the Jones polynomial also arises from the spin preserving $ U_q(sl_2) \ R$-matrix, 
Kauffman said that the Jones and Alexander polynomials have the same footing \cite{LK}. 
In this subsection, we construct a Yang-Baxter state model for $ R_{1.2} $ to show 
the relevance of the $R$-matrix for the Alexander polynomial. 
The $R_{1.1}$ case will be studied elsewhere. 
  
  A Yang-Baxter state models is a combinatorial summation that is well-defined 
on oriented diagrams. For a given diagram of a link, states are produced based on the 
$R$-matrix. Each state carries a uniquely determined quantity and their 
sum is shown to be an invariant. We mention that the first state model 
for the Alexander polynomial, that is not Yang-Baxter state model, is 
given in \cite{LK2}. The Yang-Baxter state model discussed in \cite{LK} 
is technically distinct from the ones for other invariants. We have to 
represent links as two-strand tangles, otherwise the invariant vanishes 
for all links. The reason of this is similar to Turaev's construction. 
Naive construction of the state model shows that the value for loops has 
to be zero. We therefore regularize the unknot by replacing it with 
an unknotted strand and assign $1$ for the strand.  

\begin{figure}[htbp] 
  \begin{center}
    \includegraphics[scale=0.5]{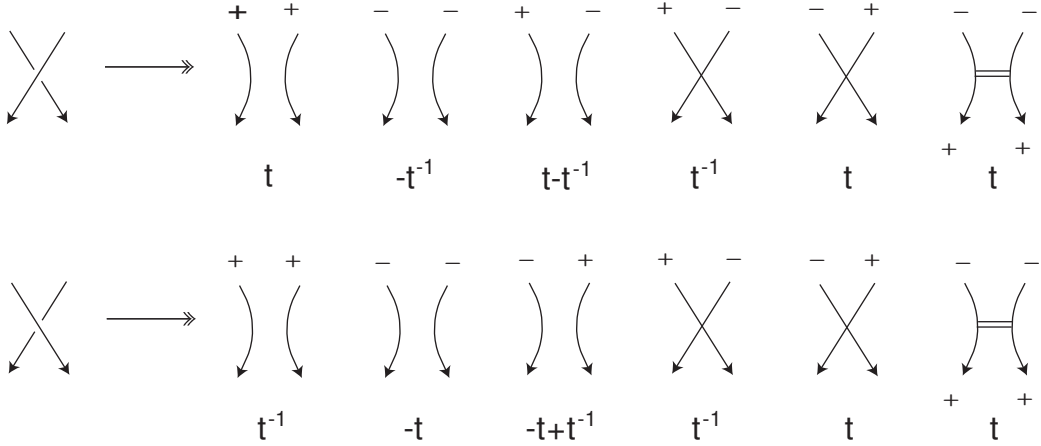}
  \end{center}
  \caption{Splicing and projection.}
  \label{crossings}
\end{figure}

  We follow this regularization and investigate the simpler case $R_{1.2}$. 
Set $ R = t R_{1.2} $ and $ t = q^{-1/2}$ 
\beq
  R = \left(
    \begin{array}{cccc}
      t  & 0 & 0 & t \\
      0  & 0 & t^{-1} & 0 \\
      0  & t & t-t^{-1} & 0 \\
      0  & 0 &   0    & -t^{-1}
    \end{array}
    \right), \label{R12t}
\eeq
then
\beq
  R - R^{-1} = (t-t^{-1}) I,  \label{RRI12}
\eeq
where $I$ denotes the identity matrix. 
The $R$-matrix and its inverse are assigned to positive and negative 
crossings of diagrams, respectively. Thus the positive and 
negative crossings have the splicings and projections given in 
Figure~\ref{crossings}.
The vertex weights for each splicing and projection are indicated below the 
diagrams. The diagrams with double-line stand for spin non-preserving 
contributions. Note that the spin change occurs only from 
negative to positive. Thus the spin changing diagrams never form 
closed loops, that is, the spin changing vertices have no contributions 
to states. Let us denote the state sum for a diagram $L$ by
\beq
  \VEV{L} = \sum_{\sigma} \braket{L}{\sigma}  i^{||\sigma||},
  \label{Lstate}
\eeq
where $ \sigma $ stands for a state and $ \braket{L}{\sigma} $ denotes 
the product of vertex weights of the state. 
$||\sigma||$ is the sum of the indices of the loops in $\sigma$ multiplied 
by the rotation number defined by
\beq
   {\rm rot}\left( \raisebox{-4mm}{\includegraphics[scale=0.5]{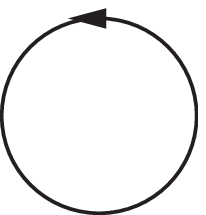}} \right) = +1,
   \qquad
   {\rm rot}\left( \raisebox{-4mm}{\includegraphics[trim=0 0 0 0,scale=0.5]{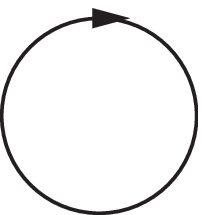}} \right) = -1.
   \label{rotnum}
\eeq
By the regularization mentioned above, for the unknotted strand
\beq
  \VEV{\includegraphics[width=15mm]{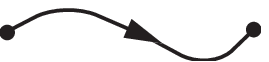}}
  = \VEV{\includegraphics[width=15mm]{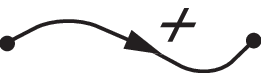}} 
  + \VEV{\includegraphics[width=15mm]{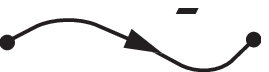}}  = 1. 
  \label{unknotted}
\eeq
It is seen by the splicings and projections given in Figure~\ref{crossings} 
that the state summation satisfies
\beq
  \VEV{ \raisebox{-5mm}{\includegraphics[width=8mm]{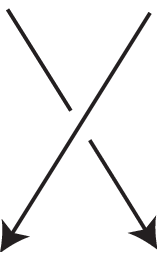}} } - 
  \VEV{ \raisebox{-5mm}{\includegraphics[width=8mm]{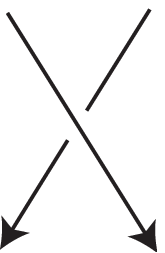}} } = 
  (t-t^{-1})  \VEV{ \raisebox{-4mm}{\includegraphics[width=8mm]{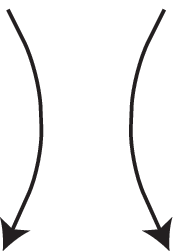}} }.
  \label{SkeinState12} 
\eeq
Since the spin non-preserving vertices do not contribute to the state summation, 
one can prove the invariance of $ \VEV{L} $ under the Reidemeister moves II, III 
essentially the same way as in \cite{LK}. It is easy to see the behaviour under 
the Reidemeister move I
\bea
  & & \VEV{ \raisebox{-2mm}{\includegraphics[width=15mm]{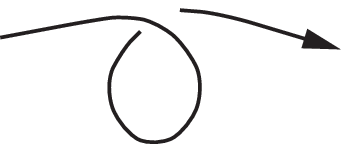}} } 
      = -i t \VEV{\includegraphics[scale=0.4, trim=0 0 0 0]{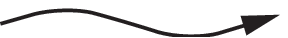}},
      \qquad
      \VEV{ \raisebox{-2mm}{\includegraphics[trim=0 0 0 0,width=15mm]{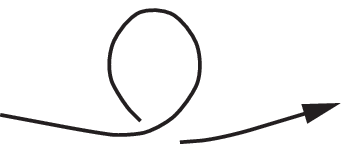}} } 
      = i t^{-1} \VEV{\includegraphics[scale=0.4, trim=0 0 0 0]{arrow.eps}},
      \nn \\
   & & \VEV{ \raisebox{-3mm}{\includegraphics[trim=0 0 0 0,width=15mm]{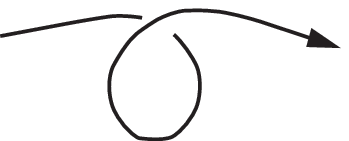}} } 
      = -i t \VEV{\includegraphics[scale=0.4, trim=0 0 0 0]{arrow.eps}},
      \qquad
      \VEV{ \raisebox{-2mm}{\includegraphics[trim=0 0 0 0,width=15mm]{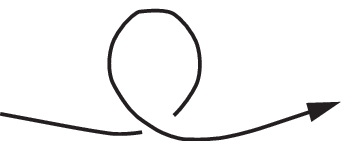}} } 
      = i t^{-1} \VEV{\includegraphics[scale=0.4, trim=0 0 0 0]{arrow.eps}}. 
      \label{RI12}
\eea 
We obtain an invariant of links by normalising $ \VEV{L} $ by 
the rotation number of the diagram $L$.
The poof of the following proposition is easy, so it is omitted. 
\begin{prop}
 \label{Alex1}{\it
   Define 
  $ \nabla_L = (it^{-1})^{-{\rm rot}(L)} \VEV{L},$ then
  \begin{enumerate}
    \item $ \nabla_L $ is an invariant of links.
    \item $ \nabla_{\includegraphics[width=8mm]{strand.eps}} = 1. $
    \item $ \nabla_{L_+} - \nabla_{L_-} = z \nabla_{L_0}, $
    where $ z = t-t^{-1}.$  
  \end{enumerate}
  }
\end{prop}
By this proposition, one can say that $ \nabla_L $ is 
indeed the Alexander polynomial. We note that $ \nabla_L $ 
is independent of the choice of the segment dropped to make 
a link diagram to a corresponding tangle. This also follows 
from the proposition \ref{Alex1} (See \cite{LK}). 

  We observed that the spin non-preserving vertex does not contribute 
to the Alexander polynomial. The legitimacy of this fact is obvious by 
comparing the spin preserving part of $ R_{1.2}$ and $ R_{2.2}.$

%
%
%
\section{Higher Dimensional Extensions}
\setcounter{equation}{0}

\subsection{Dressings of EYB}

  We learnt from the results of the last section that 
two dimensional solutions of YBE never give link invariants 
better than the Jones polynomial. Even the $R$-matrix contains 
two or three parameters, obtained invariants are reduced to 
polynomials in one variable. 
We thus have to use higher dimensional solutions of YBE to 
obtain better link invariants. However, classification of 
all solutions of YBE for $ N \geq 3 $ is far from completion. 
Even for $ N = 3, $  only the partial classification is known \cite{Hi3}. 
One way to obtain higher dimensional solutions is to use 
representations of quantum groups. Investigation of link invariants 
along the line has been done by many authors 
(see for example \cite{KR,ZGB}). 
Another way is an embedding of a known $R$-matrix to a larger matrix and then make it 
solve YBE \cite{Hi3}. Two kinds of such a method, called \textit{dressings}, are 
presented in \cite{Hi3}. In this section, we apply the Turaev's construction 
to $R$-matrices obtained by dressings and study the constructed link invariants. 
We will see that improvement of link invariants by dressings is minor. 
Thus we should find solutions of YBE not by dressings if we want better invariants. 

  Let us define two kinds of dressings. Let $ {\cal I} = \{ \ 1, 2, \cdots, N \ \} $ 
be a set of indices and $ {\cal J} $ be a selection of $ M (< N) $ numbers 
from $ {\cal I}. $ Let $ \tR $ be a $M$ dimensional solution of YBE and 
$R$ be a $N$ dimensional dressed solution.

\noindent
(a) diagonal dressing
\beq
   R^{k\ell}_{ij} = \left\{
      \begin{array}{lcl}
         \tR^{k\ell}_{ij} & \quad & i, j, k, \ell \in {\cal J} \\
         s_{ji} \delta^{k}_j \delta^{\ell}_i & & {\rm otherwise}
      \end{array}
   \right.
   \label{Ddressing}
\eeq
is a solution of YBE, if $ s_{im} $ satisfies
\bea
   & & \tR^{k\ell}_{ij} (s_{im} s_{jm} - s_{km} s_{\ell m}) = 0, \nn \\
   & & \tR^{k\ell}_{ij} (s_{mi} s_{km} - s_{jm} s_{m \ell}) = 0, \label{condDdressing} \\
   & & \tR^{k\ell}_{ij} (s_{m\ell} s_{mk} - s_{mj} s_{mi}) = 0. \nn
\eea

\noindent
(b) block dressing 
\beq
   R^{k\ell}_{ij} = \left\{
      \begin{array}{lcl}
         \tR^{k\ell}_{ij} & \quad & i, j, k, \ell \in {\cal J} \\
         \delta^k_j F^{\ell}_i & & i, \ell \in {\cal J},\ j, k \notin {\cal J} \\
         G^k_j \delta^{\ell}_i & & j, k \in {\cal J}, \ i, \ell \notin {\cal J} \\
         f_{ij} \delta^{k}_j \delta^{\ell}_i & & {\rm otherwise}
      \end{array}
   \right.
   \label{Bdressing}
\eeq
is a solution of YBE, if the matrices $ F, G $ satisfy
\bea
   & & ( F \otimes F ) \tR = \tR ( F \otimes F ), \qquad
       ( G \otimes G ) \tR = \tR ( G \otimes G ), \label{condBdressing} \\
   & & ( F \otimes 1 ) \tR ( G \otimes 1 ) = ( 1 \otimes G ) \tR ( 1 \otimes F ), 
       \qquad 
       [F, \, G] = 0. \nn
\eea
An important sub-case is $ G = F^{-1}. $ In this case, the first equation of (\ref{condBdressing}) 
is sufficient. 
We note that our definitions of dressings are slightly different form the ones in \cite{Hi3}. 
One difference comes from that we use a braid group form of YBE. Another is that we put 
additional constants $ f_{ij} $ at a diagonal part of block dressing. The role of the constants 
will be clear in consideration of EYB.

  We now turn to construction of EYB for the dressed $R$-matrices. 
Let $ \tS = (\tR, \tmu, \alpha, \beta ) $ be an EYB for $\tR. $ 
Our result for a diagonal dressing is summarised as

\begin{prop}
  Let $R$ be diagonal dressings of $\tR$. Then $ S = (R, \mu, \alpha, \beta) $ 
is an EYB for the following $ \mu$'s.
\begin{enumerate}
  \renewcommand{\labelenumi}{\rm (\arabic{enumi})}
  \item if there is no further restriction on $s_{ij}$, then
  \beq
     \mu^k_j = \left\{
        \begin{array}{lcl}
           \tmu^k_j & \quad & k, j \in {\cal J} \\
           0        & & {\rm otherwise}
        \end{array}
     \right.
     \label{muDdressing1}
  \eeq
  is the only possible $ \mu. $
  \item if $ s_{kk} = \pm \alpha $ for $ k \notin {\cal J} $ and $ \tmu $ is diagonal, 
  then 
  \beq
     \mu^k_j = \left\{
        \begin{array}{lcl}
           \tmu^k_j & \quad & k, j \in {\cal J} \\
           \pm \beta \delta^k_j   & & {\rm otherwise}
        \end{array}
     \right.
     \label{muDdressing2}
  \eeq
\end{enumerate} 
\end{prop}
While our result for block dressing is as follows.

\begin{prop}
  Let $R$ be a block dressing of $\tR$. Then $ S = (R, \mu, \alpha, \beta) $ 
is an EYB for the following $ \mu$'s.
\begin{enumerate}
  \renewcommand{\labelenumi}{\rm (\arabic{enumi})}
  \item if there is no further restriction on $ F, G, f_{ij}, $ then
  \beq
     \mu^k_j = \left\{
        \begin{array}{lcl}
           \tmu^k_j & \quad & k, j \in {\cal J} \\
           0        & & {\rm otherwise}
        \end{array}
     \right.
     \label{muBdressing1}
  \eeq
  is the only possible $ \mu. $
  \item if $ [F,\, \tmu] = [G, \, \tmu] = 0$ and $ f_{kk} = \pm \alpha $ for 
  $ k \notin {\cal J},$ then
  \beq
     \mu^k_j = \left\{
        \begin{array}{lcl}
           \tmu^k_j & \quad & k, j \in {\cal J} \\
           \pm \beta \delta^k_j   & & {\rm otherwise}
        \end{array}
     \right.
     \label{muBdressing2}
  \eeq
\end{enumerate} 
\end{prop}

\noindent
{\bf Remark.} Our main assertion is (\ref{muDdressing1}) and (\ref{muBdressing1}). 
Namely, dressed $R$-matrices allow only trivially dressed EYB's. We need further restrictions on 
dressed $R$-matrices if we want to have nontrivial dressed EYB's. (\ref{muDdressing2}) 
and (\ref{muBdressing2}) are examples of such a case and they may not be the only possibility 
of dressed EYB's.

\noindent
{\bf Proof.} 
We here give a proof of both propositions, since the way of proof is similar. 
We first note that it is easy to verify that the $ \mu$'s given in the propositions~2, 3 
satisfy the definition of EYB. Therefore, our main concern is the question whether there exist 
other possible $\mu$'s which make $S$ be an EYB. 
Let us set 
\beq
  \mu^k_j = \left\{
      \begin{array}{lcl}
        \tmu^k_j & \quad & k, j \in {\cal J} \\
        t^k_j    & & {\rm otherwise}
      \end{array}
  \right.
  \label{mugeneral}
\eeq
and study possible values of $ t^k_j. $ 
In terms of matrix entries, the conditions of EYB (\ref{EYB1}), (\ref{EYB2}) and (\ref{EYB3}) 
read
\bea
   & & \sum_{k_1, k_2} R^{k_1 k_2}_{j_1 j_2} \mu^{\ell_1}_{k_1} \mu^{\ell_2}_{k_2} 
       = \sum_{k_1, k_2}  \mu^{k_1}_{j_1} \mu^{k_2}_{j_2} R^{\ell_1 \ell_2}_{k_1 k_2},
   \label{EYB1mat} \\
   & & \sum_{k_1, k_2} R^{k_1 k_2}_{j_1 j_2} \mu^{\ell_1}_{k_1} \mu^{j_2}_{k_2} 
       = \alpha \beta\, \mu^{\ell_1}_{j_1},
   \label{EYB2mat} \\
   & & \sum_{k_1, k_2} (R^{-1})^{k_1 k_2}_{j_1 j_2} \mu^{\ell_1}_{k_1} \mu^{j_2}_{k_2} 
       = \alpha^{-1} \beta\, \mu^{\ell_1}_{j_1}.
   \label{EYB3mat}
\eea
For the sake of clarity, we use Greek letters for indices not belonging to $ {\cal J} $ 
throughout the proof. 

\medskip\noindent
(a) diagonal dressing

  Let us first study (\ref{EYB1mat}). We pick up the choices of indices by which (\ref{EYB1mat}) 
is quadratic in $ t^k_j. $ 
If $ j_1, j_2, \ell_1, \ell_2 \notin {\cal J},$ then
\[
     t^{\beta_1}_{\alpha_1}\, t^{\beta_2}_{\alpha_2}\, (s_{\alpha_1 \alpha_2} - s_{\beta_1 \beta_2} ) = 0.
\]
Since this is true for arbitrary $ s_{\alpha_1 \alpha_2},$ 
$ t^{\beta}_{\alpha} $ must be diagonal :  $ t^{\beta}_{\alpha} = t_{\alpha}\, \delta^{\beta}_{\alpha}. $ 
If $ j_2, {\ell}_2 \notin {\cal J} $ and $ j_1, \ell_1 \in {\cal J},$ then 
\[
    t_{\alpha}^{\ell}\, t_{j}^{\beta}\, ( s_{\alpha j} - s_{\ell \beta}) = 0.
\]
It follows that we need to set
\beq
    t_{\alpha}^{\ell} = 0 \quad {\rm or} \quad t_{j}^{\beta} = 0. \label{condt1}
\eeq
If $ \ell_1, j_2 \notin {\cal J} $ and $ j_1, \ell_2 \in {\cal J},$ then
\[
    t_{\alpha}^{\beta_1}\, t_j^{\beta_2}\, (s_{\alpha_2 j_1} - s_{\beta_1 \beta_2}) = 0.
\]
It follows that
\beq
   t_{\alpha} = 0 \quad {\rm or} \quad t_j^{\beta} = 0.  \label{condt2}
\eeq
If $ j_1, j_2, \ell_2 \notin {\cal J} $ and $ \ell_1 \in {\cal J},$ then
\[
   t_{\alpha_1}^{\ell}\, t_{\alpha_2}^{\beta}\, ( s_{\alpha_1 \alpha_2} - s_{\ell \beta}) = 0.
\]
Thus 
\beq
   t_{\alpha}^{\ell} = 0 \quad {\rm or} \quad t_{\alpha} = 0. \label{condt3}
\eeq
Other choices of indices not belonging to $ {\cal J}$ also give one of the 
requirements (\ref{condt1})-(\ref{condt3}). 
We therefore have three possibilities:
\beq
  {\rm (i)} \ t_{\alpha}^{\ell} = t_j^{\alpha} = 0, \qquad
  {\rm (ii)} \ t_{\alpha} = t_j^{\beta} = 0, \qquad
  {\rm (iii)} \ t_{\alpha} = t_{\beta}^j = 0.
  \label{PosD}
\eeq
Note that these are not sufficient conditions for ({\ref{EYB1mat}). 
Indeed, there are other relations which must be hold to make (\ref{EYB1mat}) true. 
Those are automatically satisfied when we take into account (\ref{EYB2mat}) and (\ref{EYB3mat}). 

\noindent
Case (i) : Consider the case of $ j_1, \ell_1 \notin {\cal J} $ for (\ref{EYB2mat}) and (\ref{EYB3mat}), 
we obtain
\[
    \delta_{\nu}^{\sigma}\, t_{\nu}\, ( (s^{\pm 1})_{\nu\nu} t_{\nu} - \alpha^{\pm 1} \beta) = 0.
\]
It follows that $ t_{\nu} = 0 $ or $ t_{\nu} = \pm \beta,\ s_{\nu\nu} = \pm \alpha. $ 
The choice of $ t_{\nu} = 0 $ corresponds to (\ref{muDdressing1}). 
While for $ t_{\nu} = \pm \beta \neq 0, $ we have to return to (\ref{EYB1mat}) and 
consider the case of $ j_2, \ell_1 \notin {\cal J} $ and $ j_1, \ell_2 \in {\cal J}.$ 
For this case, (\ref{EYB1mat}) reads
\[
  t_{\alpha} \tmu^{\ell}_j (s_{\alpha j} - s_{\alpha\ell}) = 0.
\]
Since this is true for any $ s_{\alpha j}, s_{\alpha\ell}$, $\tmu^{\ell}_j$ must be diagonal.

\medskip\noindent
Case (ii) : Consider the case of $ j_1 \notin {\cal J}, \ \ell_1 \in {\cal J} $ for 
(\ref{EYB2mat}) and (\ref{EYB3mat}), we obtain
\beq
  \sum_{j \in {\cal J}} ( s_{j\nu} \tmu^{\ell}_j - \delta_j^{\ell} \alpha \beta )\, t^j_{\nu} = 0, 
  \qquad
  \sum_{j \in {\cal J}} ( (s^{-1})_{\nu j} \tmu^{\ell}_j - \delta_j^{\ell} \alpha^{-1} \beta )\, t^j_{\nu} 
  = 0.   \label{Case2D} 
\eeq
For a fixed value of $ \nu$, these are regarded as sets of linear equations in $ t^j_{\nu}.$ 
To have nonvanishing $ t^j_{\nu},$ the coefficient matrices have to be singular. 
This requirement puts a restriction on $ s_{j\nu} $ and $ s_{\nu j}. $ 
Thus we have $t^j_{\nu} = 0,$ if there is no further restriction on $s_{ji}$. 
On the other hand, if one find $ s_{j\nu} $ and $ s_{\nu j} $ 
that makes the coefficient matrices singular, there exist other $\mu$'s which is not mentioned 
in Proposition~2. However, it may be difficult to satisfy both requirements of singular coefficient 
matrices and (\ref{condDdressing}).

\medskip\noindent
Case (iii) : By considering the case of $ \ell_1 \notin {\cal J} $ and $ j_1 \in {\cal J} $ for 
(\ref{EYB2mat}) and (\ref{EYB3mat}), we obtain the same relation as 
(\ref{Case2D}) in which $ s_{j\nu} $ and $ s_{\nu j} $ are exchanged. 
We thus follow the same discussion as Case (ii). 

We have completed the proof of Proposition~2. 

\medskip\noindent
(b) block dressing

Let us first study (\ref{EYB1mat}). If $ j_1, j_2, \ell_1, \ell_2 \notin {\cal J}, $ then
\[
   t_{\alpha_1}^{\beta_1}\, t_{\alpha_2}^{\beta_2}\, (f_{\alpha_1 \alpha_2} - f_{\beta_1 \beta_2}) = 0.
\]
Thus $ t_{\alpha}^{\beta} $ must be diagonal : 
$ t_{\alpha}^{\beta} = t_{\alpha} \delta_{\alpha}^{\beta}. $
If $ j_2, \ell_1 \notin {\cal J} $ and $ j_1, \ell_2 \in {\cal J}, $ then
\[
   t_{\alpha}^{\beta}\, [F,\, \tmu]^{\ell}_j = 0.
\]
It follows that
\beq
   t_{\alpha} = 0 \qquad {\rm or} \qquad [F,\, \tmu] = 0. \label{condt4}
\eeq
If $ j_1, \ell_2 \notin {\cal J} $ and $ j_2, \ell_1 \in {\cal J}, $ then
\[
   t_{\alpha}^{\beta}\, [G,\, \tmu]^{\ell}_j = 0.
\]
It follows that
\beq
   t_{\alpha} = 0 \qquad {\rm or} \qquad [G,\, \tmu] = 0. \label{condt5}
\eeq
Thus we have two possibilities : 
\beq
   {\rm (i)} \ t_{\alpha} =0, \qquad 
   {\rm (ii)} \ [F,\, \tmu] = [G,\, \tmu] = 0.   \label{Case2B}
\eeq
We note again that these are not sufficient conditions for (\ref{EYB1mat}). 
Others are automatically satisfied when we take into account 
(\ref{EYB2mat}) and (\ref{EYB3mat}). 

\medskip\noindent
Case (i) : Consider the case of $ j_1 \notin {\cal J} $ and $ \ell_1 \in {\cal J} $ 
for (\ref{EYB2mat}) and (\ref{EYB3mat}), then we obtain
\beq
  \sum_{j,k \in {\cal J}} (G^k_j \tmu^{\ell}_k - \delta^{\ell}_j \alpha \beta )\, t^{j}_{\nu} = 0,
  \qquad
  \sum_{j,k \in {\cal J}} ((F^{-1})^k_j \tmu^{\ell}_k - \delta^{\ell}_j \alpha^{-1} \beta )\, t^j_{\nu} 
  = 0.
  \label{condt6}
\eeq
For a fixed value of $ \nu$, these are regarded as sets of linear equations in $ t^j_{\nu}.$ 
To have nonvanishing $ t^j_{\nu},$ the coefficient matrices have to be singular. 
This requirement puts a restriction on $ F, G. $ Thus we have $ t^j_{\nu} = 0, $ if 
there is no further restriction on $ F, G. $ 
Next we consider the case of $ j_1 \in {\cal J} $ and $ \ell_1 \notin {\cal J} $ for 
(\ref{EYB2mat}) and (\ref{EYB3mat}):
\beq
   \sum_{j_2,k_1,k_2 \in {\cal J}} \tR^{k_1 k_2}_{j_1 j_2}\, \tmu^{j_2}_{k_2}\, t^{\nu}_{k_1}
   = \alpha \beta\, t^{\nu}_{j_1},
   \quad
   \sum_{j_2,k_1,k_2 \in {\cal J}} (\tR^{-1})^{k_1 k_2}_{j_1 j_2}\, t^{\nu}_{k_1}\, \tmu^{j_2}_{k_2}
   = \alpha^{-1} \beta\, t^{\nu}_{j_1}.
   \label{condt7}
\eeq
To convert the equations in (\ref{condt7}) to ones containing $ F $ or $ G,$ we go back to 
(\ref{EYB1mat}). If $ \ell_1 \notin {\cal J} $ and $ j_1, j_2, \ell_2 \in {\cal J}, $ 
then we have
\[
  \sum_{k_1,k_2 \in {\cal I}} \tR^{k_1k_2}_{j_1j_2}\, t^{\nu}_{k_1} \tmu^{\ell_2}_{k_2} = 
  \sum_{k \in {\cal J}} \tmu^k_{j_1} t^{\nu}_{j_2} F^{\ell_2}_k.
\]
Taking a sum over $ j_2 = \ell_2 $ and using the first equation of (\ref{condt7}), we obtain
\beq
  \sum_{j,k \in {\cal J}} ( F^j_k\, \tmu^k_{j_1} - \delta_{j_1}^j \alpha \beta )\, t^{\nu}_{j} = 0.
  \label{condt8}
\eeq
Similarly, by considering the case  $ \ell_2 \notin {\cal J},\ j_1, j_2, \ell_1 \in {\cal J}, $ 
and taking into account the second equation in (\ref{condt7}), we obtain 
\beq
  \sum_{j,k \in {\cal J}} ( (G^{-1})^j_k\, \tmu^k_{j_1}  - \delta_{j_1}^j \alpha^{-1} \beta )\,
  t^{\nu}_j = 0.
  \label{condt9}
\eeq
For a fixed value of $ \nu, $ (\ref{condt8}) and (\ref{condt9}) are regarded as sets of 
linear equations in $ t^{\nu}_j. $ Thus the same discussion as before concludes that 
$ t^{\nu}_j = 0 $ if there is no further restriction on $ F $ and $G.$ 
If we allow to put restrictions on $F$ and $G$, there is a possibility of the 
existence of nonvanishing $ t^j_{\nu} $ and  $ t^{\nu}_j $ which is not mentioned in 
Proposition~3. However, it may be difficult to satisfy the above requirements of singular 
coefficient matrices and (\ref{condBdressing}) simultaneously. 

\medskip\noindent
Case (ii) : Assume that $ j_2, \ell_1, \ell_2 \notin {\cal J}$ and $ j_1 \in {\cal J} $ 
for (\ref{EYB1mat}), then 
\[
   t_{\alpha_2}^{\beta_2} \,
   ( \sum_{k \in {\cal J}} F_{j}^k t_k^{\beta_1} - f_{\beta_1 \beta_2} t_j^{\beta_1} ) = 0.
\]
While by assuming that $ j_1, j_2, \ell_2 \notin {\cal J} $ and $ \ell_1 \in {\cal J},$ 
we obtain
\[
   t_{\alpha_1}^{\beta_2} \,
   ( \sum_{k \in {\cal J}} G_{k}^{\ell} t^k_{\alpha_2} - f_{\alpha_1 \alpha_2} t_{\alpha_2}^{\ell} ) = 0.
\]
Since $ t_{\alpha} \neq  0, $ the following relations are required to make (\ref{EYB1mat}) true
\beq
  \sum_{k \in {\cal J}} ( F_j^k - f_{\beta_1 \beta_2} \delta_j^k )\, t^{\beta_1}_k = 0, 
  \qquad
  \sum_{k \in {\cal J}} ( G_k^{\ell} - f_{\alpha_1 \alpha_2} \delta_k^{\ell} )\, t^{k}_{\alpha_2} = 0, 
  \label{condBf}
\eeq
We set $ t^{\beta}_k = t^k_{\alpha} = 0 $ to obtain a simple nontrivial EYB. 
For nonvanishing $ t^{\beta}_k, t^k_{\alpha}, $ 
we have to put further relations between $ F $ and $f$ to satisfy (\ref{condBf}). 
Under the condition of vanishing $ t^{\beta}_k, t^k_{\alpha}, $ 
consider the case of $ j_1, \ell_1 \notin {\cal J} $ for (\ref{EYB2mat}) and (\ref{EYB3mat})
\beq
  t^{\beta_1}_{\alpha_1}\, (f_{\alpha_1 \alpha_1} t_{\alpha_1} - \alpha \beta) = 0,
  \qquad
  t^{\beta_1}_{\alpha_1}\, (f^{-1}_{\alpha_1 \alpha_1} t_{\alpha_1} - \alpha^{-1} \beta) = 0.
  \label{condBf2}
\eeq
Since $ t_{\alpha} \neq  0, $ we obtain
\beq
   f_{\alpha_1 \alpha_1} = \pm \alpha.  \label{condBf3}
\eeq
Substituting this into (\ref{condBf2}), we obtain $ t_{\alpha} = \pm \beta. $ 

  This completes the proof of Proposition~3.

\subsection{Examples of Dressed Invariants}

   In this section, we give some examples of invariants obtained from dressings of 
two dimensional EYB's given in Table~1. 

The cases of (\ref{muDdressing1}) and (\ref{muBdressing1}) are trivial, $i.e.$ 
no new invariant is obtained. To see this, note that 
$ \displaystyle \sum_{k_1,k_2} R^{k_1 k_2}_{j_1 j_2} \mu_{k_1}^{\ell_1} \mu_{k_2}^{\ell_2} = 0 $ 
if at least one of $ \ell_1, \ell_2 $ is not in the index set $ {\cal J}.$ 
It follows that 
$ {\rm Tr}(\pi(\xi) \circ \mu^{\otimes n} ) = {\rm Tr}(\tilde{\pi}(\xi) \circ \tmu^{\otimes n}), $ 
where $ \pi(\sigma_i^{\pm 1}) = R^{\pm 1}_i $ and $ \tilde{\pi}(\sigma_i^{\pm 1}) = \tR^{\pm 1}_i. $ 
Thus $ T_S(L) = T_{\tilde{S}}(L). $ 

   We, therefore, have to consider the cases of (\ref{muDdressing2}) 
and (\ref{muBdressing2}) to obtain nontrivially dressed invariants. 
We first study  extensions to three dimensional $R$-matrices.  
Let $ {\cal I} = \{ \, 1, 2, 3 \ \} $ and $ {\cal J} = \{ \ 1, 3 \ \}. $ Throughout this section, 
we set $ \beta = 1. $ 
Note that all the EYB's given in Table~1 allows only diagonal dressings, since 
the conditions $ [F,\, \tmu ] = [G,\,\tmu] = 0 $ for $ \tmu$ in Table~1 require that 
$ F $ and $G$ are diagonal. This is a special case of diagonal dressing. 

\medskip \noindent
(i) $ \tR = R_{2,1}, \ \tmu = {\rm diag}(\sqrt{pq}, 1/\sqrt{pq}), \ \alpha = 1/\sqrt{pq}. $

\begin{table}[htbp]
\caption{3D diagonal dressing of $ R_{2.1}.$ Invariant for knots. $ t = pq. $}
\medskip
\centerline{
\begin{tabular}{|c|c|c|}
\hline 
 {} & & Jones polynomial \\
 \cline{3-3} 
 knots  & braid  & dressed invariant \\
\hline 
 {} & & 1 \\
 \cline{3-3} 
 $0_1$ & $\sigma_11$ & 1 \\
\hline 
 {} & & $ t(1+t^2-t^3) $ \\
 \cline{3-3} 
 $3_1$ & $\sigma_1^3$ & $ t^{1/2}(1-t+2t^{3/2}-t^2+t^3-t^{7/2}) $ \\
\hline 
 {} & & $ t^{-2} (1-t+t^2-t^3+t^4) $ \\
 \cline{3-3}
 $ 4_1$ & $\sigma_1 \sigma_2^{-1} \sigma_1 \sigma_2^{-1}$ & 
 $ t^{-2} (1 - t^{1/2} + t^{3/2} - t^2 + t^{5/2} - t^{7/2} + t^4) $ \\
\hline 
 {} & & $ t^{2} (1+t^2-t^3+t^4-t^5) $ \\
 \cline{3-3}  
 $ 5_1 $ & $ \sigma_1^5 $ & 
 $ t^{1/2}(1-t^{1/2}+2t^{3/2} - 2t^2 + t^{5/2} + t^3 - t^{7/2} + t^{9/2} -t^{5} + t^{6} - t^{13/2}) $ \\
\hline 
 {} & & $ t(1-t+2t^2-t^3+t^4-t^5) $ \\
 \cline{3-3} 
 $ 5_2 $ & $ \sigma_2^2\sigma_1^{-1} \sigma_2 \sigma_1^{2} $ & 
 $ t^{1/2} (1-t+t^{3/2}+t^{7/2}-t^4+t^{5} - t^{11/2}) $ \\
\hline
\end{tabular}
}
\label{R21knot}
\end{table}
This is the case of $ T_{\tS}(L) $ being Jones polynomial. 
A diagonal dressing for $\tS $ contains three additional parameters. We set 
\beq
   s_{12} = by, \quad s_{23} = ay, \quad s_{21} = a, \quad s_{32} = b, \quad 
   s_{22} = \frac{1}{\sqrt{pq}}.
   \label{DdressR12}
\eeq
The dressed invariants for some knots and links are shown in Table~2 and Table~3, 
respectively. The Jones polynomial computed from $R_{2.1}$ is also shown 
in the tables for comparison.  
It is observed that the dressing makes Jones polynomial 
more complex and the additional parameters appear only in links.

\begin{table}[htbp]
\caption{3D diagonal dressing of $R_{2.1}.$ Invariant for links. $ s = aby. $ 
The dressed invariant is not 
normalised by the value of unknot: $ t^{1/2} + 1 + t^{-1/2}.$ }
\medskip
\centerline{
\begin{tabular}{|c|c|c|}
\hline 
 {} & & Jones polynomial \\
 \cline{3-3} 
 links  & braids & dressed invariant \\
\hline 
 {} & & $t^{1/2}(1+t^2)$ \\
 \cline{3-3} 
 $2^2_1$ & $ \sigma_1^2$ & $ 2+t+t^2+t^3+2st^{1/2}(1+t)$ \\
\hline 
 {} & & $ t^{3/2}(1+t^2-t^3+t^4) $ \\
 \cline{3-3} 
 $4^2_1 $ & $ \sigma_1^4$ & $ 1+t+t^2+t^3+t^6 + 2s^2t^{3/2}(1+t) $ \\
\hline 
 {} & & $ -t^{-7/2} (1-2t+t^2-2t^3+t^4-t^5) $ \\
 \cline{3-3} 
 $ 5^2_1$ & $ \sigma_1\sigma_2^{-1} \sigma_1 \sigma_2^{-2} $ & 
 $ -t^{-4}(1-t-t^2-t^3-2t^4-t^6-2t^{7/2}-2t^{9/2})$ \\
\hline 
 {} & & $t^{5/2} (1+t^2-t^3+t^4-t^5+t^6)$ \\
 \cline{3-3} 
 $ 6^2_1 $ & $ \sigma_1^6 $ & 
 $ 1+t^2+t^3+t^4+t^9 + 2s^3 t^{5/2}(1+t) $ \\
\hline 
 {} & & $t^{3/2} (1-t+2t^2-2t^3+2t^4-t^5+t^6)$ \\
 \cline{3-3} 
 $ 6^2_2 $ & $ \sigma_2^3\sigma_1^2\sigma_2\sigma_1^{-1} $ & 
 $ 1+t+t^3+t^6+t^8+2s^3t^{5/2}(1+t)$ \\
\hline 
 {} & & $t^{-3/2} (1-2t+2t^2-2t^3+3t^4-t^5+t^6)$ \\
 \cline{3-3} 
 $ 6^2_3 $ & $ \sigma_2\sigma_1^{-1}\sigma_{2}\sigma_3^{-1}\sigma_2\sigma_1\sigma_2\sigma_3^{-1} $ & 
 $ t^{-2}(1-t+t^2+t^4+2t^5+2t^7) + 2s^2 t^{3/2}(1+t)$ \\
\hline
\end{tabular}}
\label{R21link}
\end{table}

\noindent
(ii) $ \tR = R_{2.2}, \ \tmu = {\rm diag}(1/\sqrt{pq},-1/\sqrt{pq}),\ \alpha = \sqrt{pq}.$ 

  This is the case in which all the invariants equal to zero and have the 
skein relation of Alexander polynomial. A diagonal dressing for $\tS$ contains 
three additional parameters. We set
\beq
   s_{12} = a, \qquad s_{23} = b, \qquad s_{21} = by, \qquad s_{32} = ay, \qquad 
   s_{22} = \sqrt{pq}.
   \label{DdressR22to3}
\eeq
The dressed invariants take the same value for all knots and links 
in Table~2 and Table~3: $ T_S(L) = 1. $ 

  Next we study a four dimensional extension of $R_{2.2}.$ 
Let $ {\cal I} = \{ \ 1, 2, 3, 4 \ \} $ and $ {\cal J} = \{\ 1, 3 \ \}. $ 
$ \tS $ is the same as case (ii). We consider a diagonal dressing and this 
case allows us to have ten additional parameters. We set
\bea
  & & s_{12} = a, \quad s_{14} = c, \quad s_{23} = b, \quad s_{24} = h, \quad s_{34} = d, 
  \qquad  s_{21} = by, \nn \\
  & &  s_{41} = dw, \quad s_{32} = ay, \quad s_{42} = g, \quad s_{43} = cw, \quad s_{22}=s_{44} = \sqrt{pq}. 
  \label{DdressingR22to4} 
\eea
The dressed invariants for all knots given in Table~2 take the same value: $ T_S(L) = 1. $ 
While there is a slight improvement for links shown in Table~4.

\begin{table}[htbp]
\caption{4D diagonal of $R_{2.2}.$ Invariant for links. $ t = pq, \ s = hg. $}
\medskip
\centerline{
\begin{tabular}{|c|c|c|c|c|}
\hline 
  link & $ 2^2_1$ & $4^2_1(4) $ & $4^2_1(4) $ & $ 5^2_1 $ \\
\hline  
  invariant & $ 1 + st^{-1} $ & $1 + (st^{-1})^2$ & $1 + (st^{-1})^2$ & 2 \\
\hline
\end{tabular}}
\label{R22link4D}
\end{table}

%
%
%
\section{Concluding Remarks}	
\setcounter{equation}{0}

  We have seen that two dimensional solutions of YBE produce 
no link invariants better than the Jones polynomial. An interesting 
observation is that when a solution contains two parameters, those parameters 
are combined to give one parameter polynomial invariant. 
This means that we have to use higher dimensional solutions 
to obtain more powerful polynomial invariants. One way of constructing 
higher dimensional solutions of YBE is dressings of known lower dimensional 
ones. We learnt via some examples that improvement of invariants by dressing of 
two dimensional solutions is minor. This does not mean that dressings are useless 
to find powerful invariants. Dressings of higher dimensional $R$-matrices 
have a possibility to produce better invariants. A problem of two dimensional $R$-matrices 
is that they allow only diagonal dressings. Block dressings will appear if we start 
with a higher dimensional $R$-matrix and we anticipate to obtain better invariants, 
since block dressings give more complexity to EYB. 

   Another interesting result of the present work is that there exist three 
two dimensional $R$-matrices which produce the Alexander polynomial. The spin preserving 
one ($R_{2.2}$) is a well-known $R$-matrix which has a connection to the Alexander polynomial. 
We have observed that two spin non-preserving $R$-matrices also have a connection 
to the Alexander polynomial. Spin non-preserving part of one of them ($R_{1.2}$) does not 
contribute to computation of the polynomial, while other's ($R_{1.1}$) does. 
It would be an interesting problem to construct a state model for $R_{1.1}$, since 
a role of spin non-preserving part will be clear in the computation of state models. 
It will be a future work.

%
%
%
\section*{Appendix}

 This Appendix is devoted to show the following fact.
\begin{prop}
 Let $ S = (R, \mu, \alpha, \beta) $ be an EYB and $ R' $ be a solution 
of YBE obtained from $R$ by one of the transformations (\ref{similarity}) - (\ref{LRflip}). 
Then $ S' = (R', \mu', \alpha', \beta) $ is also an EYB for the following choice of 
$ \mu', \alpha' $ and $ T_{S'}(L) = T_S(L). $
\begin{enumerate}
\renewcommand{\labelenumi}{(\arabic{enumi})}
   \item $ \mu' = \kappa Q \mu Q^{-1}, \ \alpha' = \kappa \alpha \ $ for (\ref{similarity})
   \item $ \mu'^k_i = \mu^i_k, \ \alpha' = \alpha \ $ for (\ref{transpose})
   \item $ \mu'^k_i = \mu^{k+n}_{i+n}, \ \alpha' = \alpha \ $ for (\ref{shift})
   \item $ \mu' = \mu, \ \alpha' = \alpha \ $ for (\ref{LRflip})
\end{enumerate}
\end{prop}

\noindent
{\bf Proof.} It is straightforward and easy to verify that $ S' $ of (1)-(3) satisfy the definition of 
EYB. Case (4) requires a bit care. The lhs of (\ref{EYB2}) and (\ref{EYB3}) for $S'$ of (4) 
are reduced to
\beq
   {\rm Tr}_2 (R'^{\pm 1} \circ (\mu' \otimes \mu')) = {\rm Tr}_1 (R^{\pm 1} \circ (\mu \otimes \mu)) 
   \label{Tr2Tr1}
\eeq 
It can be shown that one can change $ {\rm Tr}_2 $  to $ {\rm Tr}_1 $ in the definition of EYB and 
$ T_S(L) $ defined by (\ref{Ts}) is invariant again. We assume for case (4) that 
$ S = (R, \mu, \alpha, \beta) $ is EYB defined with $ Tr_1.$ 
Then, one can show that $S'$ of case (4) becomes EYB with $Tr_2$.

  We now turn to the proof of $ T_{S'}(L) = T_S(L). $ 
Case (3) is obvious, since $ T_S(L) $ is obtained by summing over all indices. 
Case (2) is an exchange of upper indices and lower ones. This corresponds to vertical 
flip of the braid $ \xi \in B_n $ whose closure is isotopic to $L$. Thus $ T_{S'}(\xi) $ is equal to 
$ T_S(\eta),$ where $ \eta \in B_n $ is the braid obtained by flipping $ \xi $ vertically. 
Since $ \bar{\eta} = \bar{\xi},$ we obtain $ T_{S'}(L) = T_S(L). $ 
Similarly, case (4) is corresponding to horizontal flipping of a braid. 
For case (1), it is easy to see
\[
  T_{S'}(L) = \alpha'^{-w(\xi)} \beta^{-n} \kappa^{w(\xi)} 
            {\rm Tr}(Q^{\otimes n}\, \circ \pi(\xi)\circ \mu^{\otimes n} \circ
            (Q^{-1})^{\otimes n}),
\]
where $ \pi(\sigma_i^{\pm 1}) = R_i^{\pm 1}. $ 
It follows immediately that $ T_{S'}(L) = T_S(L). $

%
%
\section*{Acknowledgements}

  One of the authors (N.A.) is partially supported by the grants-in-aid 
from JSPS, Japan (Contract No. 15540132).


\begin{thebibliography}{99}
\bibitem{Tu} V. G. Turaev, \textit{ The Yang-Baxter equation and invariants 
of links}, Invent. Math. {\bf 92} (1988) 527--553.
\bibitem{LK} L. H. Kauffman, \textit{ Knots and Physics}, World Scientific (1993).
\bibitem{AW1} Y. Akutsu and M. Wadati, \textit{ Knot invariants and the critical statistical systems}, 
J. Phys. Soc. Jpn. {\bf 56} (1987) 839.
\bibitem{AW2} Y. Akutsu and M. Wadati, \textit{ Knots, links, braids and exactly solvable models in 
statistical mechanics},  
Comm. Math. Phys. {\bf 117} (1988) 243. 
\bibitem{AW3} Y. Akutsu and M. Wadati, \textit{ Exactly solvable models and new link polynomials. I. 
N-state vertex models}, J. Phys. Soc. Jpn. {\bf 56} (1987) 3039--3051.
\bibitem{ADW1} Y. Akutsu, T. Deguchi and M. Wadati, \textit{ Exactly solvable models and new link 
polynomials. II. 
Link polynomials for closed 3-braids}, J. Phys. Soc. Jpn. {\bf 56} (1987) 3464--3479.
\bibitem{DAW} T. Deguchi, Y. Akutsu and M. Wadati, \textit{ Exactly solvable models and new link 
polynomials. III. 
Two-variable topological invariants}, J. Phys. Soc. Jpn. {\bf 57} (1988) 757--776.
\bibitem{ADW2} Y. Akutsu, T. Deguchi and M. Wadati, \textit{ Exactly solvable models and new link 
polynomials. IV. 
IRF models}, J. Phys. Soc. Jpn. {\bf 57} (1988) 1173--1185.
\bibitem{DWA} T. Deguchi, M. Wadati and Y. Akutsu, \textit{ Exactly solvable models and new link 
polynomials. V. 
Yang-Baxter operators and braid-monoid algebra}, J. Phys. Soc. Jpn. {\bf 57} (1988) 1905--1923.
\bibitem{CP} V. Chari and A. Pressley, \textit{ A Guide to Quantum Groups}, Cambridge Univ. 
Press (1994). 
\bibitem{Hi1} J. Hietarinta, \textit{ All solutions to the constant 
quantum Yang-Baxter equation in two dimensions}, Phys. Lett. {\bf A165} 
(1992) 245--251.
\bibitem{Hi2} J. Hietarinta, \textit{ Solving the two-dimensional constant 
quantum Yang-Baxter equation}, J. Math. Phys. {\bf 34} (1993) 1725--1756.
\bibitem{SUAW} K. Sogo, M. Uchinami, Y. Akutsu and M. Wadati, \textit{ Classification 
of exactly solvable two-components models}, Prog. Theor. Phys. {\bf 68} (1982) 508--526.
\bibitem{Ja} F. Jaeger, \textit{ Composition products and models for the 
HOMFLY polynomial}, L'Enseignment Math. {\bf 35} (1989) 323--361.
\bibitem{KS} L. H. Kauffman and H. Saleur, \textit{ Free fermions and the 
Alexander-Conway polynomial}, Comm. Math. Phys. {\bf 141} (1991) 293--327.
\bibitem{RS} L. Rozansky and H. Saleur, \textit{ Quantum field theory 
for the multi-variable Alexander-Conway polynomial}, Nucl. Phys. {\bf B376} 
(1992) 461--509.
\bibitem{LK2} L. H. Kauffman, \textit{ Formal Knot Theory}, Mathematical Notes 30, 
Princeton University Press (1983).
\bibitem{Hi3} J. Hietarinta, \textit{ The upper triangular solutions to the three-state 
constant quantum Yang-Baxter equation}, Preprint, solv-int/9306001.
\bibitem{KR} A. N. Kirillov and N. Yu. Reshetikhin, \textit{ 
Representations of the algebra $U_q(sl(2))$, $q$-orthogonal polynomials and invariants of links}, in 
\textit{ Infinite dimensional Lie algebras and groups}, Adv. Ser. Math. Phys. {\bf 7} (1989) 285--339.
\bibitem{ZGB} R. B. Zhang, M. D. Gould and A. J. Bracken, \textit{ Quantum group invariants 
and link polynomials}, Comm. Math. Phys. {\bf 137} (1991) 13--27.
\end{thebibliography}
\end{document}